\documentclass[11pt]{article}
\usepackage{graphicx} 
\usepackage{amsfonts} 
\usepackage{amssymb}  
\usepackage{amsmath}
\usepackage{epstopdf}
\usepackage{color}
\usepackage[utf8]{inputenc}

\newtheorem{proposition}{Proposition}
\newtheorem{lemma}{Lemma}

\begin{document}
\title{\bf{Isogeometric simulation of acoustic radiation}}

\author{Victoria Hern\'andez Mederos $^a$, Eduardo Moreno Hern\'andez $^a$, Jorge Estrada Sarlabous $^a$,\\
Isidro A. Abell\'o Ugalde $^b$, Domenico Lahaye $^c$\\
$^a$ {\it Instituto de Cibern\'etica, Matem\'atica y F\'isica, ICIMAF, La Habana, Cuba}\\
$^b$ {\it CEPES, Universidad de La Habana, Cuba}\\
$^c$ {\it DIAM, TU Delft, The Netherlands}
}
\date{}
\maketitle

\begin{abstract}
In this paper we discuss the numerical solution on a simple 2D domain of the Helmoltz equation with mixed boundary conditions. The so called {\it radiation problem} depends on the wavenumber constant parameter $k$ and it is inspired here by medical applications, where a transducer emits a pulse at a given frequency. This problem has been successfully solved in the past with the classical Finite Element Method (FEM) for relative small values of $k$. But in modern applications the values of $k$ can be of order of thousands and FEM faces up several numerical difficulties. To overcome these difficulties we solve the radiation problem using the Isogeometric Analysis (IgA), a kind of generalization of FEM. Starting with the variational formulation of the radiation problem, we show  with details how to apply the isogeometric approach in order to compute the coefficients of the approximated solution of radiation problem in terms of the B-spline basis functions. Our implementation of IgA using GeoPDEs software shows that isogeometric approach is superior than FEM, since it is able to reduce significatively the pollution error, especially for high values of $k$, producing additionally smoother solutions which depend on less degrees of freedom.
\end{abstract}


\section{Introduction}

Wave problems have been intensively studied since they are relevant in multiple fields.
The solution of wave equation is usually written as the product of a function
of time and a function $u(x,y)$ which only depends on spacial variables. In acoustic problems, for instance, the time function is chosen as  $e^{\textrm{i}\omega\,t}$, where $\omega$  is the angular frequency of the propagating wave and \textrm{i} is the imaginary unit.
With this time harmonic dependence, the wave equation is reduced to the Helmholtz equation
$$-\triangle u(x,y) - k^2u(x,y)=0 $$
where $k = 2\pi/\lambda$ is the number of waves per unit of distance, called {\it wavenumber}, and $\lambda$ is the {\it wavelength}.

Helmholtz equation is very important in acoustic applications, including seismic wave propagation, acoustic noise control, non destructive testing and medical ultrasound. In particular, therapeutic applications of ultrasound involve focused beams directing the ultrasound energy into the tissue region that needs the treatment. Currently, High Intensity Focused Ultrasound (HIFU) therapy method is known as one of the most advanced surgical \cite{Iza20} and  also physio-therapeutical techniques \cite{Gut12}. In most clinical applications, HIFU transducers are excited at a single frequency in the range $0.5-8$ MHz. From the mathematical point of view, Pennes' bioheat equation \cite{Pen48} is used to model thermal diffusion effects of HIFU. This equation relates the temperature distribution
in time and space with the absorbed ultrasound energy,
which is computed from the acoustic pressure field $u(x,y)$ solution of the Helmholtz equation.

The numerical solution of Helmholtz equation is in general a challenge. When the wavenumber $k$ is small, it can be handled using low order Finite Element Method (FEM). But the design of robust and efficient numerical algorithms for high values of $k$ is difficult. In practice, many numerical difficulties appear. First, since the function $u(x,y)$ oscillates on a scale of $1/k$, to obtain an accurate approximation of $u(x,y)$ using finite elements of degree $p$ and increasing $k$, it is necessary to request that the total number $N$ of degrees of freedom to be proportional to $k^2$ \cite{Gand15}, or to choose a mesh size $h$ such that $hk^{(p+1)/p}$ is constant and sufficiently small \cite{Mel10}. It means that very large linear systems have to be solved, with high computational cost. Moreover, it is known \cite{Bab95}, \cite{Bab00}, \cite{Ihl95}, that even if we use a big number of degrees of freedom, the errors of continuous Galerkin finite element approximations increases when $k$ becomes larger. In the literature, this non-robust behavior with respect to $k$, is known as the {\it pollution} effect.
To reduce the pollution several authors have proposed to enrich the basis of the finite element space with wave-like functions depending on the wavenumber. One of the first steps in this direction is the Partition of Unity Finite Element Method \cite{Mel96}, \cite{Moh10}.

The standard variational formulation of the Helmholtz equation is sign-indefinite (i.e.not coercive). Hence, another difficulty for the numerical solution of the Helmholtz equation is that for $k$ sufficiently large, the coefficient matrix is indefinite and non-normal. As a consequence, iterative methods to solve the corresponding linear systems behave extremely bad if the system is not preconditioned \cite{Ern11}, \cite{Diw20}. To face this problem researches have proposed several preconditioners, such as multigrid methods with Krylov smoothers,  domain decomposition, and complex shifted Laplacian preconditioner. The last one was introduced in \cite{Erl04} and further developed and  successfully generalized in
\cite{She13}, \cite{Gand15} and \cite{She16}.

Dealing with wave problems, the small discrepancies between the boundary of the mesh constructed by FEM and the boundary physical domain $\Omega$, can significantly increase the error of the FEM approximated solution \cite{Moh10}. This is more evident in 3D industrial applications, where the surface of the physical domain is usually represented in terms of Nonuniform Rational B-spline functions (NURBS) \cite{Piegl97}. Since B-spline spaces include as a particular case the piecewise polynomial spaces commonly used in FEM,
it was natural to think of the possibility of writing the approximated solution of the partial differential equation (PDE) in terms of the B-spline basis functions. This idea led to the emergence of the Isogeometric Analysis (IgA), introduced by Hughes et al.  in 2005 \cite{Hug05}, as a modern method to solve PDE. IgA uses B-spline functions to parametrize the geometry $\Omega$ {\it and} as shape functions to approximate the solution of the PDE. In this sense, it combines the variational techniques of isoparametric FEM,  with the classical functions in computer design systems. IgA and FEM are based on the same principle, the Galerkin method, but  IgA approach has a very important advantage: B-spline basis functions may be constructed to have high smoothness. This is crucial in problems with smooth solutions, where improved accuracy per degree of freedom is obtained in comparison with the classical FEM. It explains the wide range of applications solved successfully  in the last years with IgA approach, see for instance \cite{Cott06}, \cite{Baz06}, \cite{Zha07}, \cite{Buff10}, \cite{Her18}.

\subsection{Related work}
The literature dealing with different aspects of the solution of Helmholtz equation with IgA is very recent \cite{Khaj16},\cite{Coox16},\cite{Diw19}, \cite{Diw20} and \cite{Dwa21}.

In \cite{Khaj16}  the performance of IgA to solve exterior scattering problems is investigated, using an absorbing boundary condition on a fictitious boundary to truncate the infinite space. It is shown that IgA is a robust approach to reduce the effects of the pollution error and therefore it is a promising tool to solve high-frequency acoustic problems.
In \cite{Coox16} IgA is used to solve Helmholtz equation with several boundary conditions in 2D regions. The results of a convergence study are presented confirming that IgA outperforms FEM for similar degrees of freedom, specially when the frequency of the waves increases. In \cite{Diw19} the effect of higher continuity of B-spline basis function on the pollution error is studied. The conclusion is that the pollution is improved with IgA compared to classical FEM. Moreover, it is shown that partition of unity isogeometric analysis (PUIgA) is suitable for wave problems, since enrichment eliminates the need of domain re-meshing at higher frequencies.

In \cite{Diw20} the Helmholtz equation with Robin boundary condition is tackled using IgA. GMRES method for solving the linear system resulting from IgA is investigated, including the use of preconditioners such as ILU with a complex shift and complex shifted Laplace. The study in \cite{Diw20} concludes that, for all wavenumbers,  GMRES converges at a fewer iterations with IgA compared to FEM. Moreover, the pollution error is significantly reduced with IgA, even when it is not completely eliminated. In a very recent paper \cite{Dwa21}, the focus of the research is on the numerical solution of the linear system derived from IgA discretization of Helmholtz equation. The system is solved with GMRES and its convergency is accelerated using a deflation technique, combined with the approximated computation of the inverse of the CSLP with a geometric multigrid method. Numerical results for one a two dimensional problems are shown, confirming scalable convergence with respect to the wavenumber and the order of the B-spline basis functions.

\subsection{Our contribution}
The main contribution of this paper is the application of IgA to the solution of a radiation problem, mathematically modeled with a  2D Helmholtz equation with mixed boundary conditions. All the details concerning the application of IgA method are presented, including the computation of the matrix and the right hand side vector of the linear system, whose solution provides the B-spline coefficients of the approximated solution. With the open source software GeoPDEs \cite{Fal11} we have implemented our in-house code to solve the radiation problem using IgA.  The results of the numerical implementation with IgA are compared with the solution using classic FEM, for relative high values of the wavenumber. We show the superiority of the performance of IgA approach by means of several experiments, which confirm that using less degrees of freedom smoother approximated solutions are obtained with a substantially reduced pollution.


\section{Physical problem and variational formulation}

\subsection{The radiation problem}

In this paper we are interested in acoustic wave applications. Under the assumption that the acoustic wave propagation is linear and
also that the amplitude of shear waves in the media are much smaller than the amplitude of the pressure waves, nonlinear effects and
shear waves may be neglected. In consequence, the {\it acoustic wave pressure} $u(x,y)$ is a complex function
solution of the Helmholtz equation and the wavenumber $k$ is a positive and real number.

\begin{figure}[htb]
\center
\includegraphics[scale=0.35]{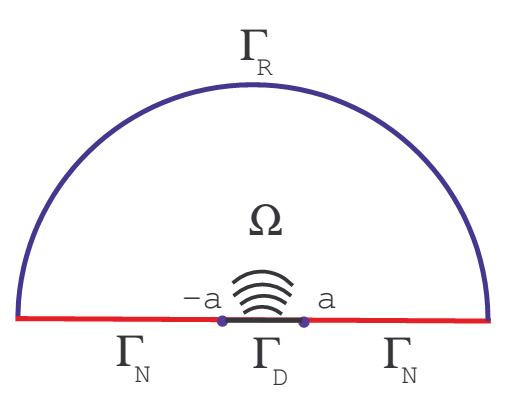}
\caption{Physical domain of radiation problem. The boundary $\Gamma$ is subdivided in three curves: $\Gamma_D, \Gamma_N$ and $\Gamma_R$. A transducer of aperture $2a$ is located in the origin of coordinates. }
\label{Fig:Omega}
\end{figure}

Inspired by the experiments to measure focused ultrasound induced heating in a tissue phantom \cite{Baz09}, \cite{Mar15}
we consider a 2D axial symmetry geometry $\Omega$: the semicircle of radius $r$ and center on the origin of coordinates.
Moreover, we assume that a transducer of aperture $2a$, with $0 < a < r$, and flat geometry is located at the origin, see Figure \ref{Fig:Omega}. The transducer emits a piston-like pulse of frequency $f=\frac{c\, k}{2\pi}$, with speed $c$ and constant amplitude
equal to $C>0$. Hence, Dirichlet boundary condition $u(x,y)=C$ is imposed on $\Gamma_D:=\{(x,y), |x| \leq a, \; y=0\}$.
Additionally, boundary $\Gamma_N:=\{(x,y), a <|x| < r, \; y=0\}$ is simulated as an acoustically rigid wall by setting the normal
velocity equal to zero. Dirichlet and Neumann boundary conditions  are known in the literature as rigid and free baffle respectively.
Finally, assuming that $a<<r$ we require that the boundary $\Gamma_R:=\{(x,y)/ x^2+y^2=r,\; y>0\}$ has the same acoustic impedance as the media  to avoid wave reflections. This Robin condition is referred to as impedance boundary condition.

In the rest of the paper we call {\it radiation problem} to the solution of equation
\begin{equation}
-\triangle u(x,y) - k^2u(x,y)=0, \;\;\;(x,y)\in \Omega \label{Heq}
\end{equation}
with {\it mixed} boundary conditions
\begin{eqnarray}
u(x,y)&=& C  \;\;\;\; \mbox{on} \,\,\Gamma_D \label{Dirichlet} \\
\frac{\partial u(x,y)}{\partial \overrightarrow{n}} &=&  0 \;\;\;\; \mbox{on} \,\,\Gamma_N \label{Neuman}\\
\frac{\partial u(x,y)}{\partial \overrightarrow{n}}+ \textrm{i}\, k u(x,y)&=&0 \;\;\;\; \mbox{on} \,\,\Gamma_R \label{Robin}
\end{eqnarray}
where $\overrightarrow{n}$ denotes the normal vector to the boundary $\Gamma_N$ or $\Gamma_R$.

Near the transducer, in the {\it near field} area, there are significant fluctuations in the ultrasound intensity. However, from some point on the pressure waves form a relatively uniform front that spreads out in a pattern originating from the center of the transducer. This area is called the {\it far field} and it is important in applications, since optimal detection occur at the start of far field, where the sound wave is well behaved and attains its maximum strength. The near field length $N_f$ defines the transition point between the near field and the far field. This point, sometimes referred to as the ``natural focus", can be calculated as
\begin{equation}
N_f=\frac{a^2}{\lambda}
\label{foco}
\end{equation}
where $\lambda=\frac{c}{f}=\frac{2\pi}{k}$ is the wavelength.

\subsection{Variational formulation}

In this section we obtain variational formulation of the radiation problem.
Denote by $V_0$ the Hilbert space
\begin{equation}
V_0=\{v:\Omega \rightarrow \mathbb{C}, v \in H^1(\Omega),\;v(x,y)=0\;\;\mbox{for}\;\;(x,y) \in \Gamma_D\} \label {VC}
\end{equation}
with the norm
\begin{equation}
\|v\|^2_{V_0}:=\|\nabla v\|_{L^2(\Omega)}^2=\int\int_{\Omega} \nabla v^t \nabla{\overline{v}} \; d\Omega \label{normV}
\end{equation}
where $\overline{v}$ is the complex conjugate of $v$ and $\nabla v =\left(\frac{\partial v}{\partial x},\frac{\partial v}{\partial y}\right)^{t}$ is the gradient vector of the complex function $v$. To obtain the variational formulation we multiply (\ref{Heq}) by $\overline{v}$ with $v \in V_0$ and integrate on $\Omega$. Applying Green's first identity and imposing the mixed boundary conditions we arrive to the {\it variational formulation},
\begin{equation}
\mbox{Find}\,\, u \in H^1(\Omega)\,\,\mbox{with}\,\, u=C \,\,\mbox{on}\,\Gamma_D \,\,\mbox{, such that}\,\,\,\, a(u,v)=0 \,\,\,\,\mbox{for all}\,\, v\in V_0
\label{variational}
\end{equation}
where the sesquilinear form $a(u,v): H^1(\Omega) \times H^1(\Omega) \rightarrow \mathbb{C}\,\,$  is given by,
\begin{eqnarray}
a(u,v)&=&\int\int_{\Omega} (\nabla u(x,y)^{t}\nabla \overline{v}(x,y)-k^2u(x,y)\overline{v}(x,y))\;d\Omega +\textrm{i}k\int_{\Gamma_R} u(x,y)\overline{v}(x,y) \;ds
\label{auv}
\end{eqnarray}

The solution $u$ of the variational problem (\ref{variational}) is called {\it weak} solution. We prove that it exits considering two auxiliary problems and writing the weak solution $u$ as the sum of the solutions of the auxiliary problems.

\begin{lemma}\label{Lema}
For each $f\in L^2(\Omega)$ there exits a weak solution of the auxiliary problem
\begin{equation}
-\triangle u_0(x,y) - k^2u_0(x,y)=f(x,y), \;\;\;(x,y)\in \Omega
\label{Heqf}
\end{equation}
with {\it homogeneous mixed} boundary conditions
\begin{equation}
u_0(x,y)= 0 \,\, \mbox{on} \;\;\Gamma_D, \hspace{1cm}
\frac{\partial u_0(x,y)}{\partial \overrightarrow{n}} =  0 \,\, \mbox{on} \;\; \Gamma_N, \hspace{1cm}
\frac{\partial u_0(x,y)}{\partial \overrightarrow{n}}+ \textrm{i} k u_0(x,y)=0 \;\; \mbox{on} \,\,\Gamma_R
\label{boundcond0}
\end{equation}
\end{lemma}

\smallskip

{\bf Proof} (main steps): A weak solution of the auxiliary problem is a solution of the variational problem
\begin{equation}
\mbox{Find}\,\, u_0 \in V_0\,\,\mbox{, such that}\,\,\,\, a(u_0,v)=\langle f,v \rangle \,\,\,\,\mbox{for all}\,\, v\in V_0
\label{variationalu0}
\end{equation}
where $\langle f,v \rangle$ is the scalar product in $L^2(\Omega)$ and $a(u,v)$ is given by (\ref{auv}). The form $a(u,v)$ in (\ref{auv}) is sesquilinear and continuous. Continuity follows from the Cauchy-Schwarz inequality and the continuity of the trace map $Tr: V_0 \rightarrow L^2(\Gamma)$, since  $a(u,v)$  involves an integral over $\Gamma_R \subset \Gamma$ with $\Gamma=\partial\Omega$.  Although $a(u,v)$ is not coercive, it satisfies a G{\aa}rding inequality . Indeed,
$$\mathfrak{Re}(a(v, v))+k^2 {\|v\|^2_{L^2(\Omega)}}=\|\nabla v\|^2_{L^2(\Omega)}=\|v\|^2_{V_0}.$$
Hence Fredholm theory can be applied, to show that a solution to the variational problem (\ref{variationalu0}) exists since the {\it homogeneous adjoint} problem has an unique solution. Details of the proof can be found in \cite{Her20}.
$\;\;\;\;\blacksquare$
\bigskip

The next result shows that a solution of (\ref{variational}) may be easily constructed in terms of a new function $u_C$ and the
weak solution $u_0$ of the problem (\ref{Heqf}) with $f(x,y)=k^2u_C$.

\begin{proposition}
Let $u_C\in H^1(\Omega)$  be a function satisfying the conditions
\begin{equation}
u_C(x,y)= C \;\; \mbox{on} \,\,\Gamma_D, \hspace{1cm} u_C(x,y)= 0 \;\; \mbox{on} \,\,\Gamma_R, \hspace{1cm}
\nabla u_C(x,y)=\mathbf{0}\;\; \mbox{on} \,\,\Omega
\label{conduC}
\end{equation}
Denote by $u_0$ the solution of (\ref{variationalu0}) with $f(x,y)=k^2u_C$. Then, the function
\begin{equation}
u=u_0+u_C
\label{u0puC}
\end{equation}
is a  weak solution of  (\ref{Heq}) with boundary conditions (\ref{Dirichlet})-(\ref{Robin}), i.e a solution of (\ref{variational}).
\end{proposition}

\medskip
{\bf Proof} \\
For $(x,y) \in \Gamma_D$ we have $u(x,y)=C$. Furthermore, for all $v \in V_0$, from (\ref{variationalu0}) we get,
$a(u_0,v)=\langle f,v \rangle=k^2\langle u_C,v \rangle$
. Moreover, from (\ref{auv}) and (\ref{conduC}) it holds,
$a(u_C,v)=-k^2 \int\int_{\Omega} u_C(x,y)\overline{v}(x,y))\;d\Omega =-k^2 \langle u_C,v \rangle$
. Finally, since  $a(u,v)$ is a sesquilinear form, from
the previous equalities  and (\ref{u0puC}) we obtain,
$a(u,v)=a(u_0,v)+a(u_C,v)
=0$
for all $v \in V_0$, i.e. $u$ is a solution of (\ref{variational}).
$\;\;\;\;\blacksquare$

\smallskip

{\it Observation}: function $u_C$ satisfying conditions (\ref{conduC}) is not unique. In consequence, the weak solution
of the radiation problem is not unique either.

\section{Galerkin method with isogeometric approach}

The Galerkin method solves the variational problem assuming that the approximated solution belongs to a finite-dimensional subspace $V^h$ of $V_0$. In the classical FEM, $V^h$ consists of piecewise polynomials functions with global $C^0$ continuity. In the isogeometric approach \cite{Cott09}, $V^h$ is generated by tensor product NURBS functions with higher global continuity. Moreover, the physical domain $\Omega$ is previously parametrized by a smooth function
$$\mathbf{F}(\xi,\eta):\hat{\Omega}\longrightarrow\Omega$$
defined on the unit square $\hat{\Omega}$  and with piecewise smooth inverse.

\subsection{Parametrization of the domain}\label{subsecDomainpar}

\begin{figure}[htb]
\center
\includegraphics[scale=0.5]{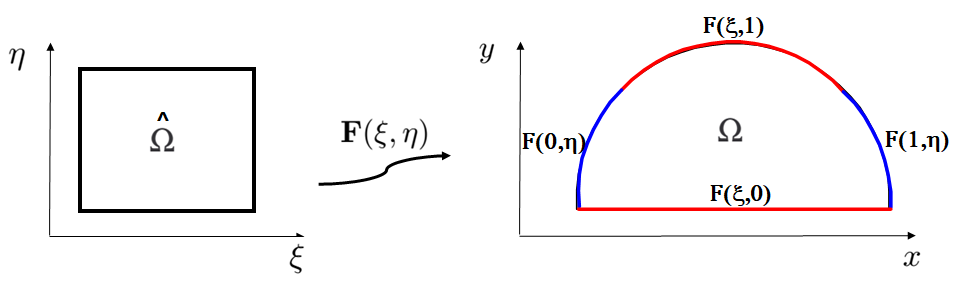}
\caption{Parametrization $\mathbf{F}(\xi,\eta)$ mapping $\hat{\Omega}$ in $\Omega$. The semicircle is subdivided in three curves.}
\label{Fig:parametrizacion3secciones}
\end{figure}

In this paper the map $\mathbf{F}(\xi,\eta)=(x(\xi,\eta),y(\xi,\eta))$ transforming $\hat{\Omega}$ in $\Omega$ is defined subdividing the semicircle in 3 circular sectors $c_l(\eta),c_r(\eta)$ and $c_t(\xi)$ with $0 \leq \xi,\eta \leq 1$, see Figure \ref{Fig:parametrizacion3secciones}. We assume that ``left" and ``right" curves $c_l(\eta)$ and $c_r(\eta)$ respectively, have both the same arc length and define an angle $\theta$ with the $x$ axis, with $0< \theta < \pi/2$. These curves can be written exactly as {\it quadratic rational  B-spline} curves for the same sequence of knots $\tau^{\eta}$. Similarly, the ``top" curve $c_t(\xi)$ can be expressed as a  quadratic rational  B-spline curve for a sequence of knots $\tau^{\xi}$.  Denote by $c_b(\xi), \; 0 \leq \xi \leq 1$ the segment of line passing through the points $(-r,0)$ and $(r,0)$. Elevating the degree of $c_b(\xi)$, it can be represented also as a  quadratic rational  B-spline curve for the sequence of knots $\tau^{\xi}$.

The map $\mathbf{F}(\xi,\eta)$ is computed as the {\it bilinearly blended Coon's patch} \cite{Piegl97} interpolating the curves $c_b(\xi)$,$c_t(\xi)$,$c_l(\eta)$ and $c_r(\eta)$. It can be written as
\begin{equation}
\mathbf{F}(\xi,\eta)=\sum_{i=1}^{n_F} \sum_{j=1}^{m_F} w_{i,j} \mathbf{P}_{i,j} R_{i,\tau^{\xi}}^{3}(\xi)R_{j,\tau^{\eta}}^3(\eta)
\label{FCoons}
\end{equation}
where $R_{i,\tau^{\xi}}^{3}(\xi),\;i=1,...,n_F$ and $R_{j,\tau^{\eta}}^3(\eta),\;j=1,..,m_F$ are respectively
the $i$-th and $j$-th rational quadratic B-splines of order $3$, for the knots $\tau^{\xi}$ and $\tau^{\eta}$, and $w_{i,j}$ are the weights.

The map $\mathbf{F}(\xi,\eta)$ satisfies the interpolation conditions
\begin{eqnarray*}
\mathbf{F}(\xi,0)&=&c_b(\xi),\;\;\mathbf{F}(\xi,1)=c_t(\xi),\;\;\;\;\;\;0 \leq \xi \leq 1 \\
\mathbf{F}(0,\eta)&=&c_l(\eta),\;\;\mathbf{F}(1,\eta)=c_r(\eta),\;\;\;\;\;\;0 \leq \eta \leq 1
\end{eqnarray*}

The points $ \mathbf{P}_{i,j},\;i=1,...,n_F,\,j=1,...,m_F$ are the vertices of the control mesh of the map $\mathbf{F}(\xi,\eta)$.
For $\theta=\pi/4$, the curves $c_l(\eta),c_r(\eta)$  can be represented as rational quadratic B-spline curves
with knots $\tau^{\eta} = \{0,0,0,1,1,1\}$, in consequence $n_F=3$. Similarly, the curve $c_t(\xi)$ can be written as a rational quadratic B-spline curve with knots $\tau^{\xi} = \{0,0,0,1,1,1\}$, hence $m_F=3$. The 9 control points $ \mathbf{P}_{i,j},\;i=1,...,3,\,j=1,...,3$
and the control mesh of $\mathbf{F}(\xi,\eta)$ are shown in second row of Figure \ref{Fig:parametrizationtheta} (left).


The quality of the parametrization $\mathbf{F}(\xi,\eta)$ of the physical domain $\Omega$ influences the precision of the solution computed with isogeometric approach \cite{Lip10}. In this sense, a good {\it uniformity} and {\it orthogonality} of the isoparametric curves of $\mathbf{F}(\xi,\eta)$ is desirable. In the literature the quality of the parametrization in the point $\mathbf{F}(\xi,\eta)$ is measured computing the {\it mean ratio Jacobian} given by,
\begin{equation}
J_r(\xi,\eta)=\frac{2\det J\mathbf{F}(\xi,\eta)}{\|\mathbf{F}_{\xi}(\xi,\eta)\|_{2}^2+\|\mathbf{F}_{\eta}(\xi,\eta)\|_{2}^2}
\label{Jacratio}
\end{equation}
where $\mathbf{F}_{\xi}(\xi,\eta)=(x_{\xi},y_{\xi})=(\frac{\partial x}{\partial \xi},\frac{\partial y}{\partial \xi})$ and
$\mathbf{F}_{\eta}(\xi,\eta)=(x_{\eta},y_{\eta})=(\frac{\partial x}{\partial \eta},\frac{\partial y}{\partial \eta})$ are the tangent vectors to the isoparametric curves, $\|\cdot\|_{2}$ denotes the Euclidean norm and $J\mathbf{F}$ denotes the Jacobian matrix of the parametrization,
$$J\mathbf{F}(\xi,\eta)=
\left(
     \begin{array}{cc}
       x_{\xi} & x_{\eta} \\
       y_{\xi} & y_{\eta} \\
     \end{array}
   \right)
$$
If the map $\mathbf{F}(\xi,\eta)$ is injective, then $\det J\mathbf{F}(\xi,\eta)$ does not changes of sign. Assuming that it is positive it holds that $0< J_r(\xi,\eta) \leq 1$. A value of $J_r$ equal to 1 at a point $\mathbf{P}_0=\mathbf{F}(\xi_0,\eta_0)$ indicates that the isoparametric curves are orthogonal at $\mathbf{P}_0$ and the map $\mathbf{F}(\xi,\eta)$  produces the same length distortion at $\mathbf{P}_0$ in both parametric directions $\xi$ and $\eta$. In the center of Figure \ref{Fig:parametrizationtheta} we show the mesh in $\Omega$ with vertices computed as the image by $\mathbf{F}(\xi,\eta)$ of the vertices of a rectangular mesh $\mathcal{Q}$ in $\hat{\Omega}$. For $\theta=\pi/4$ the corresponding parametrization $\mathbf{F}(\xi,\eta)$ preserves the geometry of the quadrilateral in $\mathcal{Q}$ is almost everywhere, except in the areas near the two points subdividing the semicircle in three curves. This is observed in the center of the second row of Figure \ref{Fig:parametrizationtheta}. On the right of this figure we show a color map, where colors correspond to the values of $J_r(\xi,\eta)$ according to (\ref{Jacratio}). A yellow color indicates that $J_r \geq 0.8$. Hence, the distortion introduced by the parametrization for $\theta=\pi/4$ is small almost everywhere, except on the two blue areas (values of $J_r$ close to $0.4$), where the distortion produced by the parametrization is higher. Observe that blue areas are not contained in the region $|x| \leq a$, where the highest pressure values are located, see Figure \ref{Fig:f1Levelcurves} (left).


\subsection{Variational formulation on $[0,1] \times [0,1]$}

With the help of the parametrization $\mathbf{F}(\xi,\eta)$ the double integral in (\ref{auv}) may be transformed into an integral over $\hat{\Omega}$ according to the integration rule
$$\int \int_{\Omega} h(x,y)\;d\Omega=\int_0^1 \int_0^1 h(\mathbf{F}(\xi,\eta))\; |\det J\mathbf{F}(\xi,\eta) |\;\;d \xi\,d \eta $$
Moreover, from the chain rule applied to $h(x,y)= h(\mathbf{F}(\xi,\eta))$ we know that $\nabla h(x, y) = J\mathbf{F}(\xi, \eta)^{-t}\nabla h(\xi,\eta)$.  Hence, the double integral in (\ref{auv}) can be written as
\begin{eqnarray}
a(u,v)&=&\int_0^1\int_0^1\left(J\mathbf{F}(\xi, \eta)^{-t}\nabla u(\mathbf{F}(\xi,\eta))\right)^{t}
\left(J\mathbf{F}(\xi, \eta)^{-t}\nabla \overline{v}(\mathbf{F}(\xi,\eta))\right) \;|\det J\mathbf{F}(\xi,\eta) |\;d \xi\,d \eta
\nonumber \\
     &-&k^2\int_0^1\int_0^1 u(\mathbf{F}(\xi, \eta))\overline{v}(\mathbf{F}(\xi, \eta))\;|\det J\mathbf{F}(\xi,\eta) |\;d \xi\,d \eta
    + \textrm{i}k\int_{\Gamma_R} u(x,y)\overline{v}(x,y) \;ds
\label{axieta1}
\end{eqnarray}

To obtain a formulation in the parametric domain it is necessary to rewrite the last integral in (\ref{axieta1}).
Recall that the parametrization $\mathbf{F}(\xi, \eta)$ has been constructed in such a way that
$\Gamma_D \cup \Gamma_N=\mathbf{F}(\xi,0),$ $\;\; 0 \leq \xi \leq 1$ and $\Gamma_R$ has been subdivided in three
consecutive boundaries $\Gamma_R=\Gamma_{R_1} \cup \Gamma_{R_2} \cup \Gamma_{R_3}$ with
\begin{eqnarray}
\Gamma_{R_1}&=&\mathbf{F}(0,\eta)=c_l(\eta),\;\; 0 \leq \eta \leq 1 \label{GR1}\\
\Gamma_{R_2}&=&\mathbf{F}(\xi,1)=c_t(\xi),\;\; 0 \leq \xi \leq 1 \label{GR2} \\
\Gamma_{R_3}&=&\mathbf{F}(1,\eta)=c_r(\eta),\;\; 0 \leq \eta \leq 1 \label{GR3}
\end{eqnarray}

Taking into account that $\mathbf{F}(\xi,\eta)=(x(\xi,\eta),y(\xi,\eta))$, from (\ref{GR1}),(\ref{GR2}) and (\ref{GR3})
we obtain an expression for $\int_{\Gamma_R} u(x,y)\overline{v}(x,y) \;ds$ that only depends on variables $\xi$ and $\eta$.
Substituting this expression in (\ref{axieta1}) we obtain finally that $a(u,v)$ can be written in the parametric domain $[0,1] \times [0,1]$ as,

\begin{eqnarray}
a(u,v)&=&\int_0^1\int_0^1 \nabla u(\mathbf{F}(\xi,\eta))^{t}(J\mathbf{F}(\xi,\eta)^{t}J\mathbf{F}(\xi,\eta))^{-1} \nabla \overline{v}(\mathbf{F}(\xi,\eta)) \;|\det J\mathbf{F}(\xi,\eta)|\;d \xi d \eta \nonumber\\
&-& k^2 \int_0^1\int_0^1 u(\mathbf{F}(\xi,\eta)) \overline{v}(\mathbf{F}(\xi,\eta))\;|\det J\mathbf{F}(\xi,\eta)|\;d \xi d \eta \nonumber\\
&+& \textrm{i}k \int_0^1\,u(\mathbf{F}(0,\eta))\overline{v}(\mathbf{F}(0,\eta))\left(\left(\frac{\partial x(0,\eta)}{\partial \eta}\right)^2 + \left(\frac{\partial y(0,\eta)}{\partial \eta}\right)^2 \right)^{1/2} \;d\eta \nonumber\\
&+& \textrm{i}k \int_0^1 \,u(\mathbf{F}(\xi,1))\overline{v}(\mathbf{F}(\xi,1))\left(\left(\frac{\partial x(\xi,1)}{\partial \xi}\right)^2 + \left(\frac{\partial y(\xi,1)}{\partial \xi}\right)^2 \right)^{1/2} \;d\xi \nonumber\\
&+& \textrm{i}k \int_0^1\,u(\mathbf{F}(1,\eta))\overline{v}(\mathbf{F}(1,\eta))\left(\left(\frac{\partial x(1,\eta)}{\partial \eta}\right)^2 + \left(\frac{\partial y(1,\eta)}{\partial \eta}\right)^2 \right)^{1/2} \;d\eta \nonumber\\
\label{axieta01R}
\end{eqnarray}

\subsection{Galerkin method with B-spline functions }

In the isogeometric approach the Galerkin method replaces the infinite dimensional space $V_0$ by the space $V^h$ generated by {\it tensor product B-spline functions} \cite{Piegl97}. These functions are computed using univariate B-splines in both parametric directions $\xi$ and $\eta$. To define the B-splines of order $k_1$ (degree $ \leq k_1-1$) in the direction $\xi$ we need a sequence of knots $t^{\xi}$. Similarly, the definition of the univariate B-splines of order $k_2$ (degree $\leq k_2-1$) in the direction $\eta$ requires a sequence of knots $t^{\eta}$. Given $n$ and $m$, the knots sequences are defined in this work by
\begin{eqnarray}
t^{\xi}&=&(\overbrace{0,...,0}^{k_1-1},\xi_1,\xi_2,...,\xi_{n-k_1+2},\overbrace{1,...,1}^{k_1-1}) \label{tchi_cuad} \\
t^{\eta}&=&(\overbrace{0,...,0}^{k_2-1},\eta_1,\eta_2,...,\eta_{m-k_2+2},\overbrace{1,...,1}^{k_2-1})\label{teta_cuad}
\end{eqnarray}
where  $0=\xi_1<\xi_2<...<\xi_{n-k_1+2}=1$ are the breakpoints in the direction $\xi$ and $0=\eta_1<\eta_2<...<\eta_{m-k_2+2}=1$ are the breakpoints in the direction $\eta$. The sequences of breakpoints define a rectangular mesh in $[0,1] \times [0,1]$ with vertices $(\xi_i,\eta_j),$ $i=1,...,n-k_1+2,\;$ $j=1,...,m-k_2+2$.

The B-splines functions $B_{i,t^{\xi}}^{k_1}(\xi)$, $i=1,...,n$ are a basis of the space $\mathbb{S}_{k_1,t^{\xi}}$ of splines of order $k_1$ with knots $t^{\xi}$. Similarly, the B-spline functions $B_{j,t^{\eta}}^{k_2}(\eta),$ $j=1,...,m$ are a basis of the space $\mathbb{S}_{k_2,t^{\eta}}$ of spline functions of order $k_2$  with knots $t^{\eta}$.  The functions,
\begin{equation}
B_{i,j}^{k_1,k_2}(\xi,\eta):=B_{i,t^{\xi}}^{k_1}(\xi)B_{j,t^{\eta}}^{k_2}(\eta),\;\;i=1,...n,\;\;j=1,...,m
\end{equation}
are a basis of the {\it tensor product space} $\mathbb{S}_{k_1,t^\xi} \bigotimes \mathbb{S}_{k_2,t^\eta}$ of splines functions of order $k_1$ in the direction $\xi$ and order $k_2$ in the direction $\eta$. To simplify the notation, in the rest of the paper we don't write the subindex $t^\xi$ or $t^\eta$ of the B-spline functions when it is clear from the context. Due to the assumptions on the parameterization $\mathbf{F}(\xi,\eta)$, the functions
\begin{equation}
\phi_{i,j}(x,y)=(B_{i,j}^{k_1,k_2}\, \circ\, \mathbf{F}^{-1})(x,y),\;\;i=1,...n,\;\;j=1,...,m
\label{Phi}
\end{equation}
are linearly independent in $\Omega$. Denote by $V^h$ the space,
\begin{equation}
V^h=span\{\phi_{i,j}(x,y),\;\;i=1,...,n,\,j=1,...,m\}
\label{Vh}
\end{equation}

Galerkin method computes the approximated solution $u^h(x,y)$ of the variational problem (\ref{variational}) as a function in $V^h$. It means that $u^h(x,y)$ is written as,
\begin{equation}
u^h(x,y)=\sum_{i=1}^{n}\sum_{j=1}^{m}\gamma_{i,j}\phi_{i,j}(x,y)
\label{uh}
\end{equation}
for certain {\it unknown} coefficients $\gamma_{i,j} \in  \mathbb{C}$. In order to obtain a system of linear equations for the unknowns $\gamma_{i,j}$ it is convenient to vectorize the basis functions and the corresponding coefficients in (\ref{uh}) introducing the
change of indices
\begin{equation}
q=q(i,j):=n(j-1)+i,\;\;i=1,...,n,\;\;j=1,...,m
\label{index}
\end{equation}
With this transformation we rewrite the expression (\ref{uh}) as
\begin{equation}
u^h(x,y)=\sum_{q=1}^{N}\alpha_q \psi_q(x,y)
\label{up}
\end{equation}
where $N:=nm$ and
\begin{eqnarray}
\mathbf{\alpha}:&=&(\alpha_1,...,\alpha_N)=(\gamma_{1,1},...,\gamma_{n,1},\gamma_{1,2},...,\gamma_{n,2},
\gamma_{1,m},...,\gamma_{n,m}) \label{alpha}\\
\mathbf{\psi}(x,y):&=&(\psi_1(x,y),...,\psi_N(x,y)) \nonumber\\
 &=&(\phi_{1,1}(x,y),...,\phi_{n,1}(x,y),\phi_{1,2}(x,y),...,\phi_{n,2}(x,y),\phi_{1,m}(x,y),...,\phi_{n,m}(x,y))
\end{eqnarray}
Thus, in the new notation
$$V^h=span\{\psi_q(x,y),\;q=1,...,N\}$$

\smallskip

Denote by $V_0^h$ the subspace of $V^h$ constituted by functions in $V^h$ {\it vanishing} on $\Gamma_D$,
\begin{equation}
V_0^h=span\{\phi_{i,j}(x,y),\;\mbox{ such that}\;\; \phi_{i,j}(x,y)=0 \;\;\mbox{for all} \;\; (x,y) \in \Gamma_D\}
\label{V0h}
\end{equation}

To determine which are the functions $\phi_{i,j}$ in $V_0^h$ we have to compute the preimage by
$\mathbf{F}(\xi,\eta)$ of the points $(-a,0)$ and $(a,0)$ delimiting $\Gamma_D$.
Since $\mathbf{F}(\xi,0)$ with $0 \leq \xi \leq 1$ is $\Gamma_D \cup \Gamma_N$, it is clear that there exist
$\xi_{a^-}$ and $\xi_{a^+}$ both in $(0,1)$ such that,
\begin{eqnarray}
\mathbf{F}(\xi_{a^-},0)&=&(-a,0) \label{chiam}\\
\mathbf{F}(\xi_{a^+},0)&=&(a,0)  \label{chiap}
\end{eqnarray}

From the expression (\ref{FCoons}) it is easy to check that the solution of equations (\ref{chiam}) and (\ref{chiap}) are:
$\xi_{a^-}=(r-a)/(2r)$ and $\xi_{a^+}=(r+a)/(2r)$, both in $(0,1)$. Assume that $\xi_{a^-}$ and $\xi_{a^+}$ satisfy

$$t_{i_1}^{\xi} \leq \xi_{a^-} < t_{i_1+1}^{\xi}, \;\;\;\;\;\; t_{i_2}^{\xi} \leq \xi_{a^+} < t_{i_2+1}^{\xi}$$

with $i_1 \leq i_2$. Since $\Gamma_D=\mathbf{F}(\xi,0)$ with $\xi_{a^-} \leq \xi \leq \xi_{a^+}$ and
$B_j^{k_2}(0)=0$, for $j=2,...,m$, we conclude that the B-splines {\it not identically null} on $\Gamma_D$ are
$$\phi_{i,1}(x,y)=B_{i,1}^{k_1,k_2}(\mathbf{F}^{-1}(x,y))=B_i^{k_1}(\xi)B_1^{k_2}(\eta),\;\;\;i_1-k_1+1 \leq i \leq i_2$$
Hence,
\begin{equation}
V_0^h=span\{\phi_{i,j}(x,y),\;\;1 \leq i \leq n,\;2 \leq j \leq m,\;\; \mbox{and} \;\; \phi_{i,1}(x,y), \;\; i \notin [i_1-k_1+1,i_2] \}
\label{V0hPhi}
\end{equation}

Denote by $I_0$ the set containing the global indices, computed using (\ref{index}), of functions on $V_0^h$. Then
$\psi_q(x,y)= 0$  for $(x,y) \in \Gamma_D$ and $q \in I_0$. In other words,
$V_0^h=span\{\psi_q(x,y),\;q \in I_0\}$ and $dim(V_0^h)=n_0$, where $n_0$ is the size of $I_0$. Similarly,
denote by $I_g$ the set of global indices of functions $\phi_{i,1}(x,y)$ with $i_1-k_1+1 \leq i \leq i_2$ and let $n_g$ be
the size of $I_g$. With this notation $n_0+n_g=N$ and $u^h(x,y)$ can be written as
$$u^h(x,y)=u_0^h(x,y)+u_g^h(x,y)$$
where $u_0^h(x,y)=\sum_{q\in I_0}\alpha_q \psi_q(x,y)$ and $u_g^h(x,y)=\sum_{q\in I_g}\alpha_q \psi_q(x,y)$. Observe that function
$u_0^h(x,y) \in V_0^h$. Hence to obtain the Galerkin formulation we substitute in (\ref{axieta01R}) $u$ by $u^h(x,y)$ given by (\ref{up}) and $v$ by a basis function $\psi_p(x,y),$ $p \in I_0$ of $V_0^h$. The result is,

\begin{equation}
\sum_{q=1}^N \alpha_q \Big( s_{p,q} \,-\,k^2\,m_{p,q} \,+\,\textrm{i}\,k\,e_{p,q} \Big) =  0\,\;\;\;\;\;\;p \in I_0
\label{eqI0}
\end{equation}
where ( omitting the dependence of functions of $(\xi,\eta)$ when it is clear)
\begin{eqnarray*}
s_{p,q}&=&\int_0^1\int_0^1 ( \nabla \psi_q)^{t}(J\mathbf{F}^{t}J\mathbf{F})^{-1} \nabla \psi_p \,|\det J\mathbf{F}|\;d \xi d \eta \\
m_{p,q}&=& \int_0^1\int_0^1 \psi_q\psi_p\,|\det J\mathbf{F}|\;d \xi d \eta \\
e_{p,q}&=& \int_0^1\,\psi_q(\mathbf{F}(0,\eta))\psi_p(\mathbf{F}(0,\eta))\left(\left(\frac{\partial x(0,\eta)}{\partial \eta}\right)^2 + \left(\frac{\partial y(0,\eta)}{\partial \eta}\right)^2 \right)^{1/2} \;d\,\eta \\
&+& \int_0^1 \,\psi_q(\mathbf{F}(\xi,1))\psi_p(\mathbf{F}(\xi,1))\left(\left(\frac{\partial x(\xi,1)}{\partial \xi}\right)^2 + \left(\frac{\partial y(\xi,1)}{\partial \xi}\right)^2 \right)^{1/2} \;d\xi \\
&+& \int_0^1\,\psi_q(\mathbf{F}(1,\eta))\psi_p(\mathbf{F}(1,\eta))\left(\left(\frac{\partial x(1,\eta)}{\partial \eta}\right)^2 + \left(\frac{\partial y(1,\eta)}{\partial \eta}\right)^2 \right)^{1/2} \;d\eta
\end{eqnarray*}
Finally we have to impose the Dirichlet boundary condition on $\Gamma_D$. For that we set,
\begin{equation}
\alpha_q=C,\;\;\;\;\;q\in I_g
\label{diag}
\end{equation}
Observe that with this assignment, for any point $(\widetilde{x},0)$ in $\Gamma_D$ it holds
$$u^h(\widetilde{x},0)=u_0^h(\widetilde{x},0)+u_g^h(\widetilde{x},0)=0+\sum_{q\in I_g} \alpha_q\psi_q(\widetilde{x},0)=C\sum_{q\in I_g}\psi_q(\widetilde{x},0)=C$$
where the last equality holds since B-spline functions satisfy the unit partition property.

Without loss of generality, assume that unknowns $\alpha_q,\;q=1,...,N$ have been reorganized in such a way that the first $n_0$ unknowns correspond to indexes in $I_0$ and the last $n_g$ unknowns correspond to indexes in $I_g$. Then, taking into account (\ref{diag}) the linear equations (\ref{eqI0}) can be written as,
\begin{equation}
\sum_{q=1}^{n_0} \alpha_q\,a_{p,q}= C\;\sum_{q=n_0+1}^{N} \,a_{p,q}   \;\;\;\;\;\;\;\;\;\;p=1,...,n_0
\label{linsystem}
\end{equation}
where $a_{p,q}=s_{p,q}\,-\, k^2\,m_{p,q}\, +\,\textrm{i}\,k\, e_{p,q}$. In the literature matrices $S=(s_{p,q})$ and $M=(m_{p,q})$ for $p,q=1,...,n_0$ are known as {\it stiffness matrix} and {\it mass matrix} respectively. Let $\widetilde{\alpha}=(\alpha_1,...,\alpha_{n_0})$, then system (\ref{linsystem}) can be written in matrix form as,
\begin{equation}
A\,\widetilde{\alpha}=b
\label{sist}
\end{equation}
where $b=(b_p)= C\;\sum_{q=n_0+1}^{N} \,a_{p,q}$ for $p=1,...,n_0 $ and $A=(a_{p,q}),\;\;p,q=1,...,n_0$ is given by
\begin{equation}
A=S\,-\,k^2M\,+\,\textrm{i}\,k\,E
\label{matrixA}
\end{equation}
with $E=(e_{p,q}),\,p,q=1,...,n_0$.

\medskip

This is a good place to recall that the numerical solution of the linear system (\ref{sist}) is a challenge. Some observations in this sense are the following. Even when matrix $A$ is symmetric, it is also indefinite and not Hermitian. Moreover, $A$ is sparse but it gets denser with the increase of the order of B-splines, since the ratio of the nonzero elements of $A$ to total number of elements of $A$ is bounded above by $(2p+1)^2/N$, with $p=\max\{k_1-1,k_2-1\}$. For large values of $k$, $A$ is also very large, owing to the fact that many degrees of freedom are necessary to obtain accurate approximations of $u(x,y)$. Therefore, direct solvers are prohibitively expensive and it is necessary to appeal to iterative solvers.

Since the convergence rate of iterative methods strongly depends on the condition number $\kappa(A)$ of matrix $A$, it is important to obtain bounds for $\kappa(A)$ as a function of the mesh size $h$ and also as function of the order of B-splines. In \cite{Gah12} bounds for the condition number of the stiffness and mass matrices of IgA discretizations of elliptic PDE are obtained. In general, it is shown that
their condition numbers grow quickly with the inverse of mesh size $h$ and the polynomial degree $p$.


As a consequence, the matrix $A$ in (\ref{matrixA}) is in general {\it ill conditioned} and the convergence of a Krylov subspace method, like  GMRES, requires a previous preconditioning of $A$. In this paper we use one of the most successful preconditioners introduced in \cite{Erl04} in FEM context, the {\it Complex Shifted Laplacian} (CSLP) $\widetilde{A}$, which depends on a parameter $\beta>0$ and it is given by
\begin{equation}
\widetilde{A}=A-\,\textrm{i}\beta M
\label{LaplacianPre}
\end{equation}
In order to obtain a good preconditioner, $\beta/k$ needs to be sufficiently small \cite{Gand15}. In \cite{Diw20} experiments are done with  $\beta=\sqrt{k}$. Moreover, to obtain wavenumber independent convergence, the complex shift $\beta$ used in \cite{Dwa21} is of order $\mathcal{O}(k)^{-1}$.

\section{Numerical results and discussion}

To solve the radiation problem, we have implemented IgA approach in an in-house code, using the open source package {\bf GeoPDEs} \cite{Fal11} to compute the matrix and the right hand side of the linear system (\ref{sist}) . IgA results reported here were obtained in a PC with i7 processor and 8Gb of RAM. In our simulations of acoustic radiation, the Dirichlet constant in (\ref{Dirichlet}) is $C=1$, the aperture of the transducer is $2a=0.02m$ and the sound propagation speed is $c=1500m/s$. We compute $\lambda=\frac{c}{f},\,$ and  $k=2\pi/\lambda$.

\subsection{Results}

In the next simulation we set the radius of the semicircle as $r=2N_f$, with $N_f$ given by (\ref{foco}). We assume that the transducer emits a piston-like pulse of frequency $f=1.0$ MHz, hence $k=4.189 \times 10^3$ and the radius of the semicircle is $r=0.133\,m$. The radiation problem is solved with {\it bicubic B-splines}, i.e $k_1=k_2=4$. Moreover, we choose $n=700$ and $m=600$ and {\it uniform} knot sequences $t^{\xi}$ and $t^{\eta}$ in the directions $\xi$ and $\eta$ respectively. Hence, the total number of degrees of freedom (dof) is  $N=n\cdot m=420\,000$. The parametrization $\mathbf{F}(\xi,\eta)$ was computed as described in section \ref{subsecDomainpar} for the parameter $\theta=\pi/4$.

GMRES is used to solve the preconditioned system
\begin{equation}
\widetilde{A}^{-1}A \widetilde{\alpha} =\widetilde{A}^{-1}b
\label{preconsys}
\end{equation}
where $\widetilde{A}$ is the CSLP (\ref{LaplacianPre}) with $\beta=1/(3k)=7.95 \times 10^{-5}$, as suggested in \cite{Dwa21}. After only 1 outer iteration of GMRES the {\it relative residual} of preconditioned system (\ref{preconsys})
was $3.729 \times 10^{-10}$.

\begin{figure}[hbt]
\center
\includegraphics[scale=0.5]{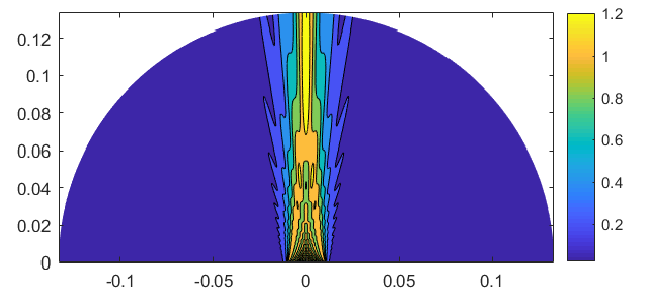}
\includegraphics[scale=0.25]{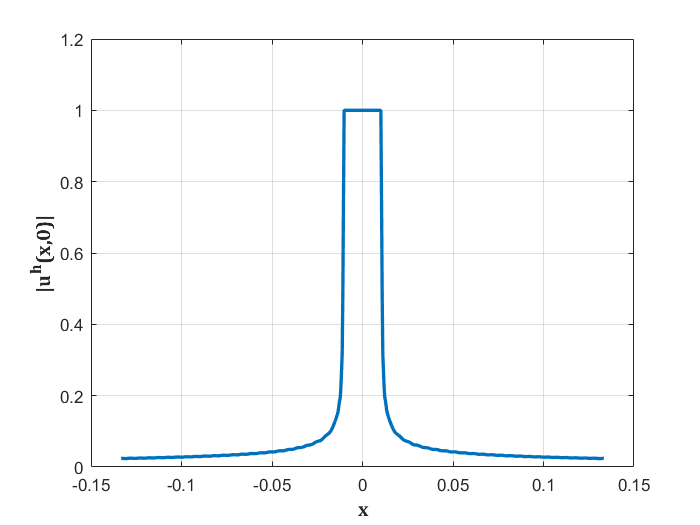}
\caption{Function $|u^h(x,y)|$ for $f=1.0$ MHz computed with bicubic B-splines with uniform knots and $N=420\,000$ dof ($n=700,\;m=600$). Right: Graphics of the function $|u^h(x,0)|$.}
\label{Fig:f1Levelcurves}
 \end{figure}

In Figure \ref{Fig:f1Levelcurves} left, we show a 2D view of the absolute value of the approximated solution,
$$|u^h(x,y)|=(({\mathfrak Re}\,u^h(x,y))^2+({\mathfrak Im}\,u^h(x,y))^2)^{\frac{1}{2}}$$ 
Colors in this figure indicate that most oscillations of acoustic pressure are in $-a \leq x \leq a$. Moreover, the region of the highest acoustic pressure has an elliptical shape and it is located after the natural focus $N_f=0.066$. In Figure \ref{Fig:f1Levelcurves} right we show the function $|u^h(x,0)|$ for $-r \leq x \leq r$. Observe that Dirichlet boundary condition holds in the interval $[-a,a]$, while in the rest of the line $y=0$ the function $|u^h(x,0)|$ is smooth.

\begin{figure}[htb]
\center
\includegraphics[scale=0.28]{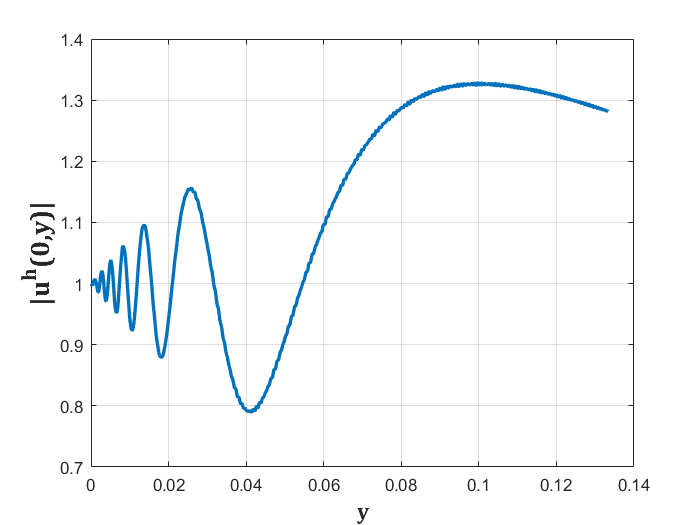}
\includegraphics[scale=0.28]{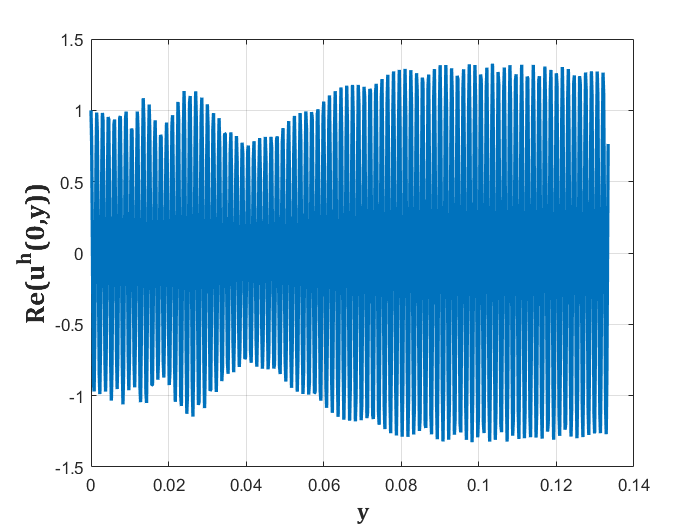}
\includegraphics[scale=0.28]{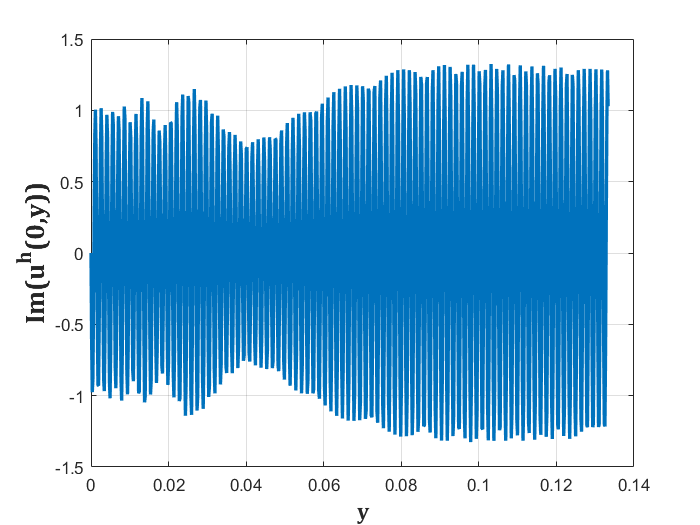}
\caption{Graphics of the functions $|u^h(0,y)|$, ${\mathfrak Re} u^h(0,y)$ and ${\mathfrak Im} u^h(0,y)$ for $f=1.0$ MHz, computed with bicubic B-splines with uniform knots and $N=420\,000$ dof ($n=700,\;m=600$).}
\label{Fig:f1cubicEjey}
 \end{figure}

It is well known that the maximum amplitude of the acoustic wave pressure is attained on the profile $x=0$. The behavior of the function $u^h(0,y)$ is shown in Figure \ref{Fig:f1cubicEjey}, where we show the graphics of the functions: ${\mathfrak Re}(u^h(0,y))$, ${\mathfrak Im}(u^h(0,y))$ and $|u^h(0,y)|$.

\subsection{Comparison with FEM}

To evaluate the performance of isogeometric approach we want to compare the IgA approximation with the solution of the radiation problem obtained using classic FEM. Unfortunately, it is difficult to compare with FEM approximations of similar problems obtained by other authors, since most of them do not report all the information needed to set up a fair comparison, such as the type of elements and the number of degrees of freedom. Hence, in this section we solve the radiation problem with FEM and IgA approaches, using in both cases piecewise {\it cubic} polynomials. FEM approximation was computed with COMSOL \cite{COM55}.  

In the next experiments we consider three values of frequency: $f=0.75$ MHz, $f=1.0$ MHz and $f=1.25$ MHz. The radius $r$ of the semicircle depends on the frequency and it is computed as $r=2N_f$. Since the wavelength $\lambda$ and the frequency $f$ are inversely proportional, if $f$ is of order $1\,MHz$, then $\lambda$ is of order $10^{-3}\,m$, which means that the solution $u(x,y)$ is highly oscillatory. To face this problem classic FEM uses a mesh of size  $h$, with $h \leq \frac{\lambda}{10}$, resulting in a high number of degrees of freedom. Here we show that the same problem can be successfully solved with the isogeometric approach with a substantially smaller number of degrees of freedom. IgA approximation is obtained using {\it cubic B-splines} with uniform knots and the parametrization $\mathbf{F}(\xi,\eta)$ described in section \ref{subsecDomainpar} for the parameter $\theta=\pi/4$. FEM approximation is computed with {\it cubic Lagrange polynomials} defined on a quadrilateral mesh. In both approaches, the linear system obtained after the discretization is solved with GMRES with Complex Shifted Laplacian Preconditioner, with parameter $\beta=1/(3k)$.

\begin{figure}[hbt]
\center
\begin{tabular}{c c c}
  $f=0.75$ MHz & $f=1.0$ MHz & $f=1.25$ MHz\\
    \includegraphics[scale=0.28]{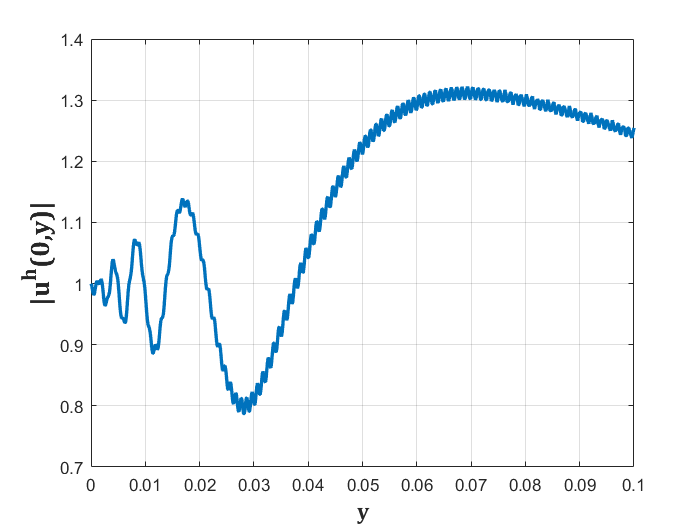}
  & \includegraphics[scale=0.28]{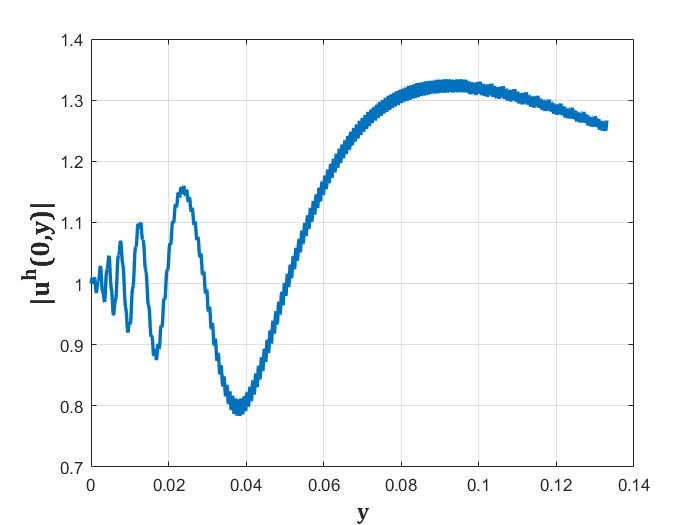}
  & \includegraphics[scale=0.28]{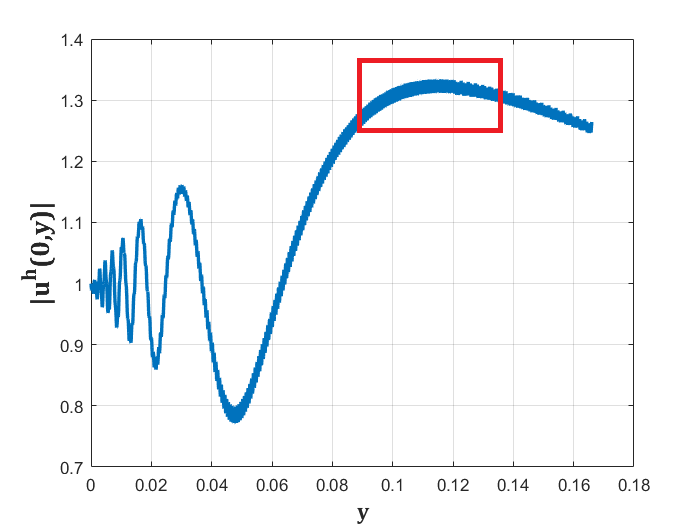}\\
  $N=1\,057\,612$ & $N=3\,330\,886$ & $N=8\,100\,010$ \\ \\
  \includegraphics[scale=0.28]{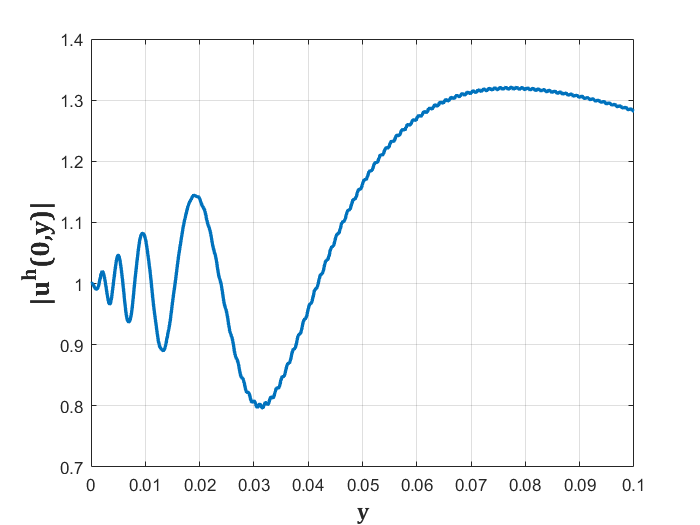} &
  \includegraphics[scale=0.28]{absuhejey1f420000dofbicubpiover4.png} &
  \includegraphics[scale=0.28]{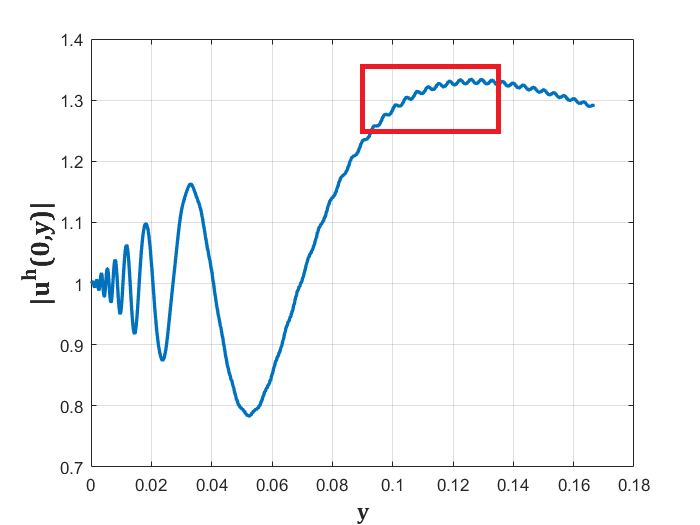}\\
  $N=200\,000$ & $N=420\,000$ & $N=720\,000$
\end{tabular}
\caption{Graphic of the function $|u^h(0,y)|$ for several values of the frequency $f$.
Top: $u^h(x,y)$ computed with cubic Lagrange FEM with a quadrilateral mesh. Bottom: $u^h(x,y)$ computed with cubic B-splines with uniform knots.  The total number of degrees of freedom $N$ is given for each case. Left: $f=0.75$ MHz, center: $f=1.0$ MHz, right: $f=1.25$ MHz. }
\label{Fig:f1compara}
 \end{figure}


\begin{figure}[hbt]
\center
\includegraphics[scale=0.28]{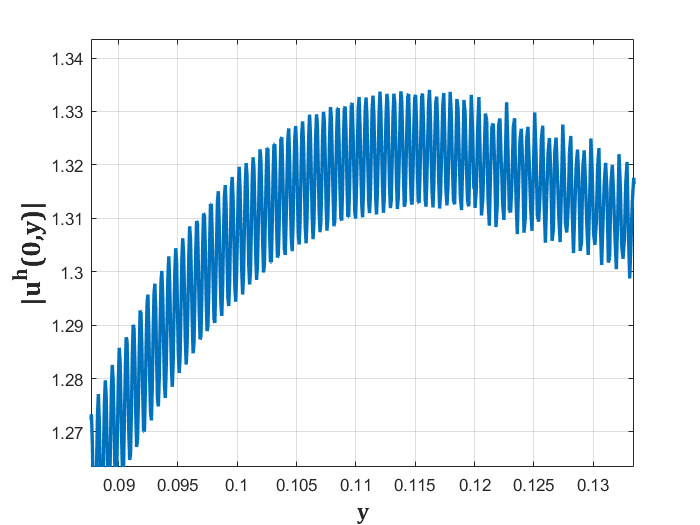}
\includegraphics[scale=0.28]{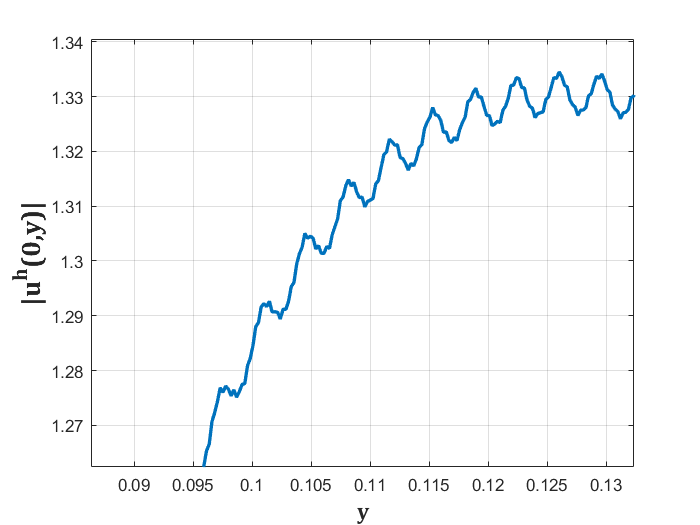}\\
\caption{Zoom of rectangular area showed in right column of Figure \ref{Fig:f1compara} ($f=1.25\,MHz$ ). Left: function $|u^h(0,y)|$ computed with FEM, right: function $|u^h(0,y)|$ computed with IgA.} 
\label{Fig:zoomf1compara}
\end{figure}

As we already mentioned, the maximum amplitude of the acoustic wave pressure is attained on the profile $x=0$, hence our comparison focus on the behavior of the function $|u^h(0,y)|$.  In Figure \ref{Fig:f1compara} we show the graphic of the functions $|u^h(0,y)|$ obtained using FEM (first row) and IgA (second row) approaches. From the physical point of view, we observe that both provide solutions with similar behavior: for increasing frequency the number of oscillations of $|u^h(0,y)|$ in the near field region grows, but in any case from the natural focus on there are no more oscillations and the maximum amplitude of the acoustic wave is reached. 

Nevertheless, comparing the functions $|u^h(0,y)|$ computed with FEM and IgA it is clear that, for all values of frequencies, FEM approximation has {\it more noise}, even when the number $N$ of degrees of freedom used for its computation is several times bigger ($5.28$, $7.93$ and $11.25$ for frequencies $f=0.75$ MHz, $f=1.0$ MHz and $f=1.25$ MHz respectively), than the number of degrees of freedom used for the computation with IgA approach. In Figure \ref{Fig:zoomf1compara} we show a zoom of the rectangular area showed in right column of Figure \ref{Fig:f1compara}. Observe that the noise of the function $|u^h(0,y)|$ computed with FEM has higher frequency and larger amplitude than the noise of the function $|u^h(0,y)|$ computed with IgA for the same frequency $f=1.25\,MHz$.

\subsection{Discussion about the IgA approach}

The practical implementation of IgA approach requires a careful selection of the degree and knots of B-spline functions in each parametric direction, the number of degrees of freedom, the parametrization $\mathbf{F}(\xi,\eta)$ of the physical domain $\Omega$  and the numerical solver of linear equations derived from the discretization. All these aspects influence the quality and accuracy of the numerical solution $u^h(x,y)$ of the PDE. In this section we show the effect of some of them.

To study the influence of the parametrization in the quality of the IgA solution we construct two parametrizations  $\mathbf{F}_{\theta_1}(\xi,\eta)$ and $\mathbf{F}_{\theta_2}(\xi,\eta)$ of the circle with radius $r=0.166$, using the method described in section \ref{subsecDomainpar} for the parameter $\theta_1=\pi/20$ and $\theta_2=\pi/4$ respectively. For the parametrization $\mathbf{F}_{\theta_1}(\xi,\eta)$, the top curve $c_t(\xi)$ belongs to the space of quadratic rational B-splines curves with knots $\tau^{\xi}=\{0, 0, 0, 1/2, 1/2, 1, 1, 1\}$. This space has dimension $n_F=5$. On the other hand, the top curve $c_t(\xi)$ for the parametrization $\mathbf{F}_{\theta_2}(\xi,\eta)$ with $\theta=\pi/4$, belongs to the space of dimension $n_F=3$  of quadratic rational B-splines curves with knots $\tau^{\xi} = \{0,0,0,1,1,1\}$.

\begin{figure}[hbt]
\center
\begin{tabular}{c c c}
\includegraphics[scale=0.28]{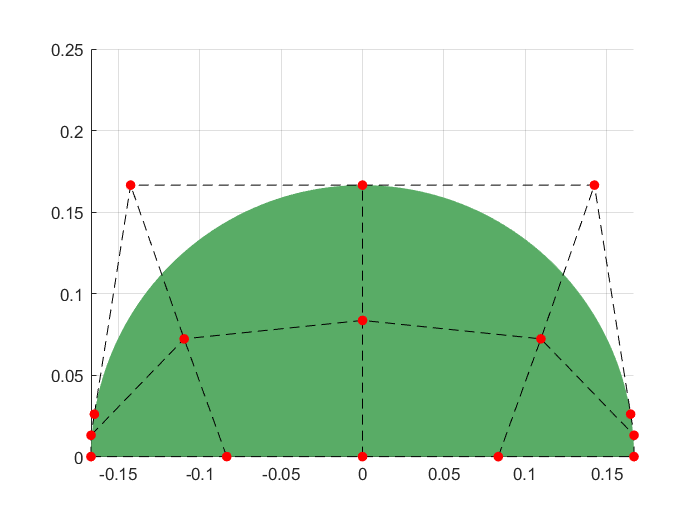} &
\includegraphics[scale=0.27]{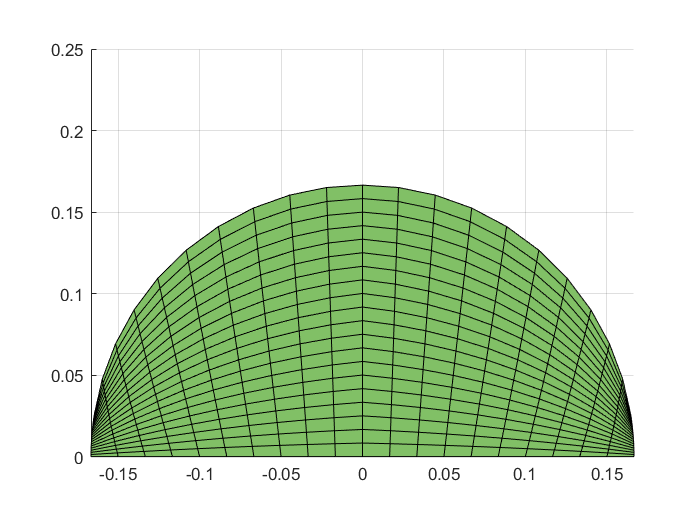} &
\includegraphics[scale=0.28]{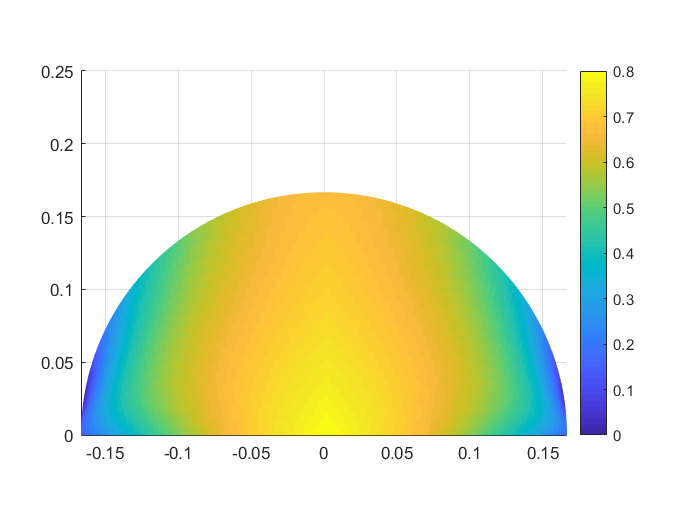}\\
\includegraphics[scale=0.28]{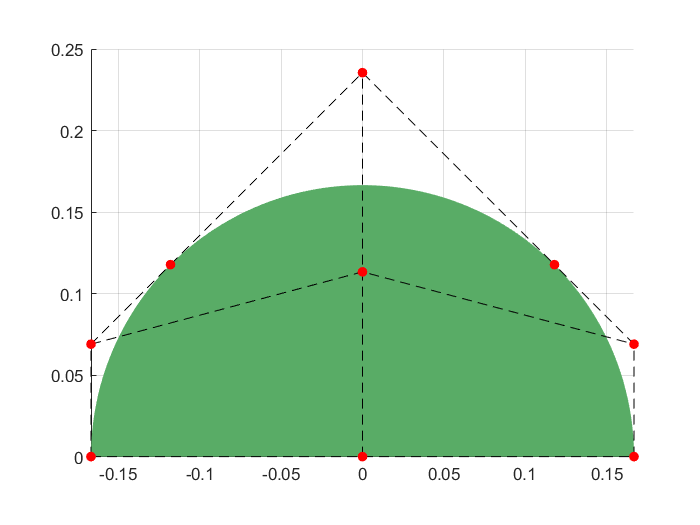} &
\includegraphics[scale=0.27]{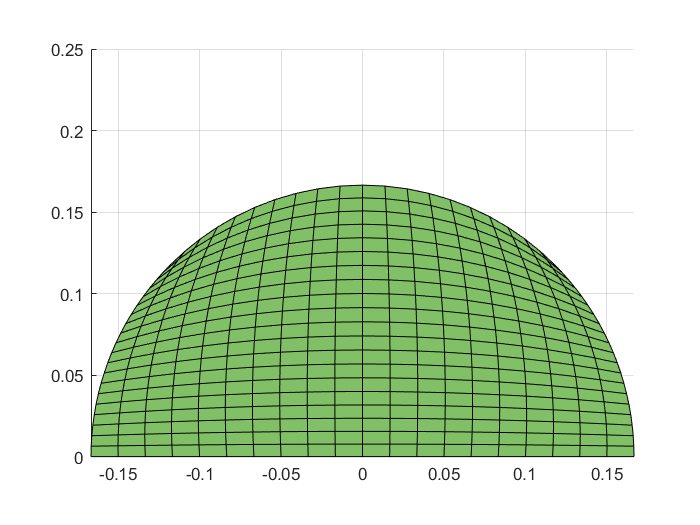} &
\includegraphics[scale=0.28]{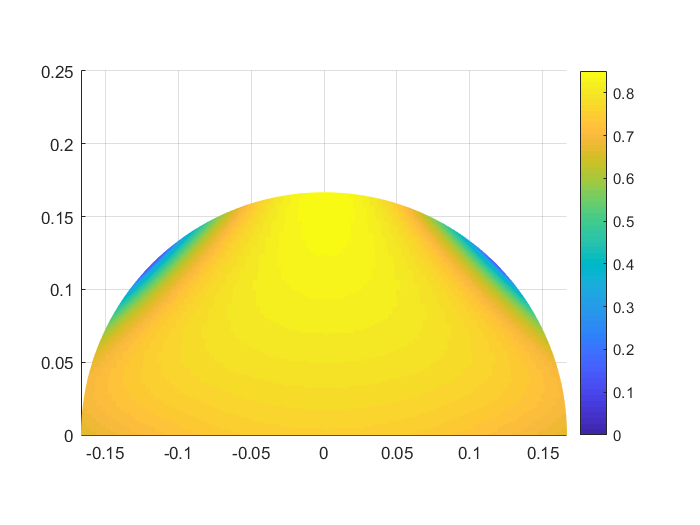}
\end{tabular}
\caption{Left: control mesh of the parametrization $\mathbf{F}_{\theta}(\xi,\eta)$, center: quadrilateral mesh obtained as image by $\mathbf{F}_{\theta}(\xi,\eta)$ of the rectangular mesh $\mathcal{Q}$ in $\hat{\Omega}$, right: mean ratio Jacobian $J_r(\xi,\eta)$. 
First row correspond to parametrization $\mathbf{F}_{\theta_1}(\xi,\eta)$ and second row to parametrization $\mathbf{F}_{\theta_2}(\xi,\eta)$.}
\label{Fig:parametrizationtheta}
\end{figure}

In each row of Figure \ref{Fig:parametrizationtheta} we show the control mesh of a parametrization $\mathbf{F}_{\theta}(\xi,\eta)$ of $\Omega$, the mesh in $\Omega$ with vertices computed as the image by $\mathbf{F}(\xi,\eta)$ of the vertices of a rectangular mesh $\mathcal{Q}$ in $\hat{\Omega}$, and a color map where colors correspond to the values of the mean ratio Jacobian $J_r(\xi,\eta)$ given by (\ref{Jacratio}) (see section \ref{subsecDomainpar}). The first row of Figure \ref{Fig:parametrizationtheta} corresponds to the parametrization $\mathbf{F}_{\theta_1}(\xi,\eta)$ and the second row to parametrization $\mathbf{F}_{\theta_2}(\xi,\eta)$. Observe that $\mathbf{F}_{\theta_2}(\xi,\eta)$ preserves the geometry of the quadrilateral mesh $\mathcal{Q}$ almost everywhere, except in the areas near the two points subdividing the semicircle in three curves. On the other hand, $\mathbf{F}_{\theta_1}(\xi,\eta)$ introduces big deformations in two areas close to the boundary curves $c_l(\eta)$ and $c_r(\eta)$. The quality of previous parametrizations is also evaluated in the color map of the mean ratio Jacobian $J_r(\xi,\eta)$. Yellow areas in this color map are those where the parametrization introduces small distortions (values of $J_r \geq 0.8$). It is clear that in this sense the parametrization $\mathbf{F}_{\theta_2}(\xi,\eta)$ is better than $\mathbf{F}_{\theta_1}(\xi,\eta)$.

To examine the effect of the parametrization and the degree of B-splines in the numerical solution $u^h(x,y)$ of the radiation problem, we solve it for $f=1.0$ MHz, first with {\it biquadratic} and later with {\it bicubic} B-splines. In both cases, we consider the previously computed  parametrizations $\mathbf{F}_{\theta_1}(\xi,\eta)$ and $\mathbf{F}_{\theta_2}(\xi,\eta)$ of $\Omega$ and $N=420\,000$ ($n=700$ and $m=600$) degrees of freedom.

\begin{figure}[hbt]
\center
\begin{tabular}{c c c}
\includegraphics[scale=0.27]{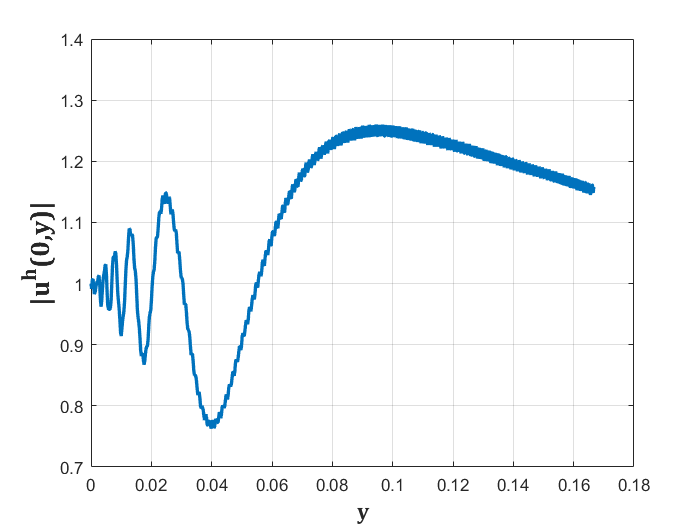} &
\includegraphics[scale=0.27]{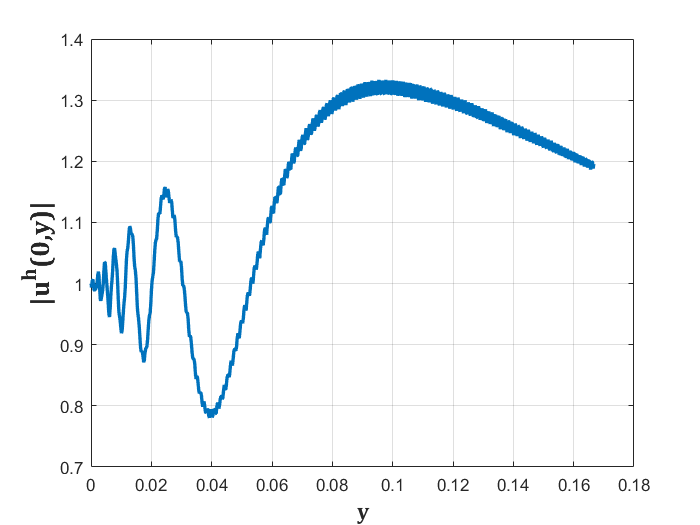}&
\includegraphics[scale=0.27]{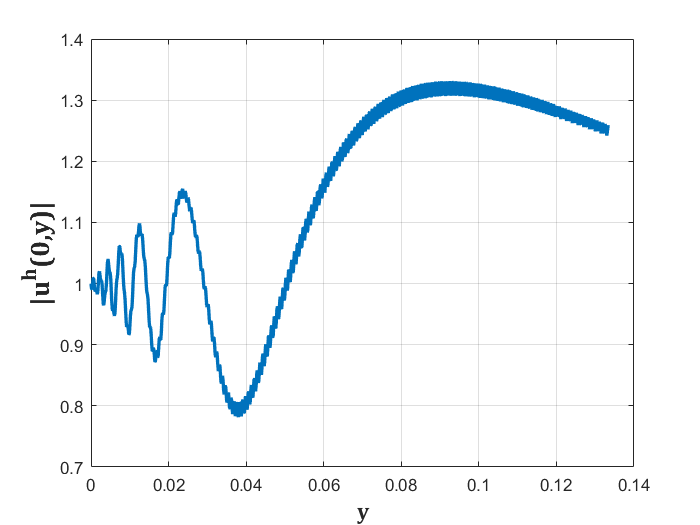}\\
quadratic  IgA, $\mathbf{F}_{\theta_1}$ & quadratic IgA, $\mathbf{F}_{\theta_2}$ & quadratic FEM \\
$N=420\,000$ & $N=420\,000$ & $N=9\,211\,077$ \\
\includegraphics[scale=0.27]{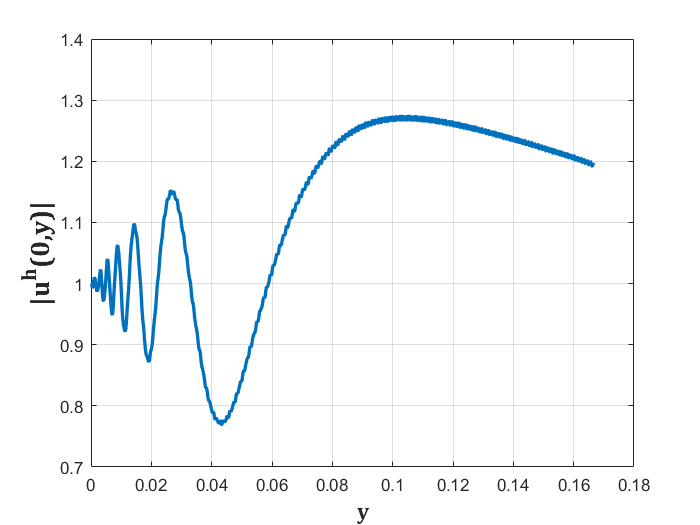}&
\includegraphics[scale=0.27]{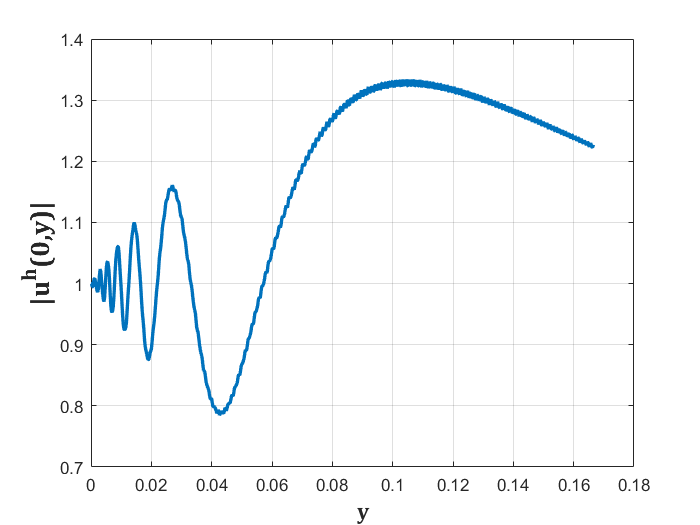} &
\includegraphics[scale=0.27]{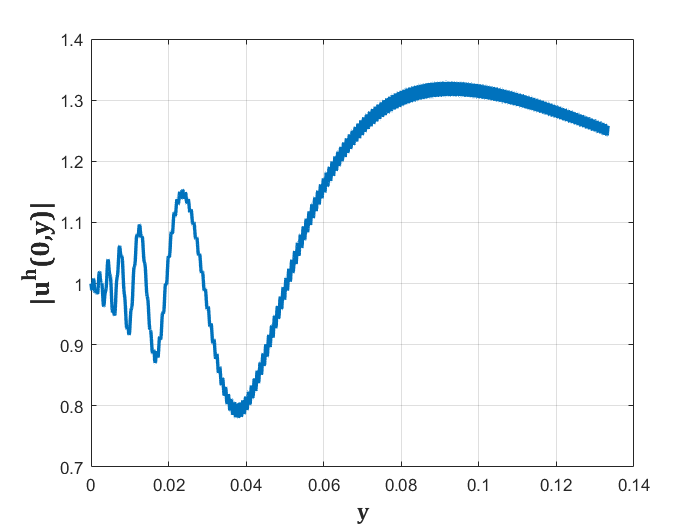}\\
 cubic IgA, $\mathbf{F}_{\theta_1}$ & cubic IgA, $\mathbf{F}_{\theta_2}$ & cubic FEM \\
 $N=420\,000$ & $N=420\,000$ & $N=20\,720\,653$
\end{tabular}
\caption{Graphic of the function $|u^h(0,y)|$ for frequency $f=1.0$ MHz.
Top: $u^h(x,y)$ computed with {\it quadratic} B-splines with uniform knots (left parametrization $\mathbf{F}_{\theta_1}$, center parametrization $\mathbf{F}_{\theta_2}$) and with {\it quadratic} FEM (right). Bottom: $u^h(x,y)$ computed with {\it cubic} B-splines  with uniform knots (left parametrization $\mathbf{F}_{\theta_1}$, center parametrization $\mathbf{F}_{\theta_2}$) and with {\it cubic} FEM (right). The total number of degrees of freedom $N$ is given for each case.}
\label{Fig:Efectothetagrado}
 \end{figure}

First row of Figure \ref{Fig:Efectothetagrado} (left and center) shows the graphic of the functions $|u^h(0,y)|$ computed with {\it quadratic} B-spline functions with uniform knots. For comparison we also show on the right the function $|u^h(0,y)|$ computed with FEM using {\it quadratic} Lagrange polynomials. Similarly, the graphic of the functions $|u^h(0,y)|$ of the second row of Figure \ref{Fig:Efectothetagrado} (left and center) were obtained computing $u^h(x,y)$ with {\it cubic} B-splines with uniform knots and, for comparison, we show on the right the function $|u^h(0,y)|$ computed with FEM using {\it cubic} Lagrange polynomials. In both rows, left and center columns correspond to IgA solutions obtained with parametrizations $\mathbf{F}_{\theta_1}(\xi,\eta)$ and $\mathbf{F}_{\theta_2}(\xi,\eta)$ respectively. Observe that as the theory predicts, for the same number of degrees of freedom, the cubic approximation is smoother than the quadratic, with less noise.

Comparing the functions $|u^h(0,y)|$ in the left and center columns of Figure \ref{Fig:Efectothetagrado}, with the function $|u^h(0,y)|$ computed with COMSOL and shown on right, we observe that IgA functions $|u^h(0,y)|$ obtained with the parametrization $\mathbf{F}_{\theta_1}(\xi,\eta)$ are poor approximations. In fact, the maximum of IgA functions $|u^h(0,y)|$ for parametrization  $\mathbf{F}_{\theta_1}(\xi,\eta)$ is 1.26 for quadratic B-splines and 1.27 for cubic, while the maximum of FEM function $|u^h(0,y)|$ is 1.33, independently of the degree of Lagrange polynomials. On the other hand, if we use the better parametrization $\mathbf{F}_{\theta_2}(\xi,\eta)$, then the IgA functions $|u^h(0,y)|$, shown in center column of Figure \ref{Fig:Efectothetagrado}, are good approximations of the function $|u^h(0,y)|$ computed with COMSOL and shown in the right column.

From similar experiments with parametrizations computed with several values of $\theta$, we conclude that the quality of the parametrization has a strong influence in the precision of the approximated solution computed with IgA. In this sense, for the parametrization $\mathbf{F}_{\theta}(\xi,\eta)$ constructed in section \ref{subsecDomainpar}, $\theta=\pi/4$ is a good option. Concerning the degree of B-splines we have observed that {\it cubic} B-splines are recommended, since with relative few degrees of freedom they allow to compute accurate IgA approximations that suffer from less noise.

\section{Conclusions}

Radiation problem appears in several important applications, where high values of the wavenumber are handled. We have solved this
problem in a simple 2D scenario, with one transducer emitting a piston-like pulse of constant amplitude. The unknown of the radiation problem is the acoustic pressure field, that from the mathematical point of view is the solution of Helmholtz equation with mixed boundary conditions.

For high values of the wavenumber, FEM solutions of radiation problem suffer from the pollution error. To overcome these limitations we have solved the radiation problem using the isogeometric method and approximating the solution with B-spline functions. Our experiments with high wavenumber have confirmed that compared to FEM, IgA shows smaller pollution errors for substantially less degrees of freedom.
The numerical solution of the linear system derived from IgA discretization is a challenge, since the matrix is large, indefinite and for large values of $k$ it is also ill conditioned. We have solved the linear system using GMRES with complex shifted Laplacian preconditioner to speed up the convergence.

In the near future we are planning to improve our implementation of the method. To overcome the restrictive tensor product structure of classic B-spline functions, adaptive isogeometric methods, like hierarchical B-splines will be used. Concerning the solution of the linear system, we will enhance the performance computing approximately the inverse of the complex shifted Laplacian preconditioner using a
geometric multigrid method. Moreover, we will tackle  more complex radiation problems, where an array of transducers emits the pulse which travels in an heterogeneous medium with wavenumber depending on the space.


\end{document}